\documentclass[10pt]{amsart}
\usepackage{amssymb,amsmath}
\usepackage[usenames]{color}

\newtheorem{thm}{Theorem}[section]

\newtheorem{lem}[thm]{Lemma}

\theoremstyle{definition}

\theoremstyle{remark}

\newcommand{\ha}{{\textstyle\frac{1}{2}}}

\def\cC{\mathcal C}
\def\cD{\mathcal D}
\def\cN{\mathcal N}
\def\cX{\mathcal X}
\def\cY{\mathcal Y}

\begin{document}

\title[System of Equations and Configurations in the Euclidean space]{System of Equations and Configurations in the Euclidean space}

\author{Annachiara Korchmaros}

\address{
Department of Computer Science \& Interdisciplinary Center for Bioinformatics\\ Universit{\"a}t Leipzig\\
H{\"a}rtelstra{\ss}e 16-18\\
D-04107 Leipzig\\
Germany}

\email{annachiara.korchmaros@uni-leipzig.de}

\date{}

\begin{abstract} 
In the $3$-dimensional Euclidean space ${\bf{E}}^3$, fix six pairwise distinct points
\begin{equation*}
\label{eqA}
\begin{array}{ccc}
A=(a_1,a_2,a_3), & B=(b_1,b_2,b_3), & C=(c_1,c_2,c_3), \\
D=(d_1,d_2,d_3), & E=(e_1,e_2,e_3), & F=(f_1,f_2,f_3)
\end{array}
\end{equation*}
together with two further points $X^*=(x_1^*,x_2^*,x_3^*)$ and $Y^*=(y_1^*,y_2^*,y_3^*)$ in $\mathbf{E}^3$. We show that System $(*)$ consisting of the following six equations
in the unknowns $X=(x_1,x_2,x_3)$ and $Y=(y_1,y_2,y_3)$
\begin{equation}
\label{egy}
\frac{1}{\|X-T\|^2}  +\frac{1}{\|Y-T\|^2}=\frac{1}{\|X^*-T\|^2} +\frac{1}{\|Y^*-T\|^2}, \quad T\in\{A,B,C,D,E,F\}
\end{equation}
has only finitely many solutions provided that both of the following two conditions are satisfied:
\begin{itemize}
\item[(i)] no four of the fixed points $A,B,C,D,E,F$ are coplanar;
\item[(ii)] no four of the six spheres of center $T$ and radius $1/\sqrt{k_T}$ with
\begin{equation}
\label{kxy}
k_T=\frac{1}{\|X^*-T\|^2} +\frac{1}{\|Y^*-T\|^2}
\end{equation}
share a common point in ${\bf{E}}^3$.
\end{itemize}
Furthermore, we exhibit configurations $ABCDEFX^*Y^*$, showing that (i) is also necessary.  
This result is an improvement on~\cite[Theorem 1]{cifuentes20233d} where the finiteness of solutions of System $(*)$ was only ensured for sufficiently generic choices of the points $A,B,\ldots,F,X^*,Y^*.$  
The extended System $(**)$ associated to System $(*)$ consists of seven equations (\ref{egy}) where $T\in\{A,B,C,D,E,E,F,G\}$ with a further point  $G=(g_1,g_2,g_3)\in{\bf{E}}^3$.
We show that if (i) and (ii) hold for $T\in\{A,B,C,D,E,F\}$ and the associated extended System $(**)$  has some solutions other
than $(X^*,Y^*)$ and $(Y^*,X^*)$, then $G$ lies on a real affine surface only depending on $\{A,B,\ldots,F\}$. This result proves \cite[Conjecture 1]{cifuentes20233d}. 

Motivation for studying the above problems comes from applications to genetics; see~\cite{cifuentes20233d}.

\end{abstract}

\maketitle

\section{Introduction}
Euclidean distance geometry is a classical topic with various interesting practical applications; see~\cite{Liberti2017}.
A problem quite often raised is determining whether a system consisting of equations depending on distances has finitely many (real) solutions. System $(*)$ has at least two (trivial)  solutions, namely $(X^*,Y^*)$ and $(Y^*,X^*)$, and the question arises whether some more solutions might exist.

System $(*)$ can have infinitely many solutions for some configurations of points $A,B,C,D,E,F, X^*$ and $Y^*$. The first problem is to address whether a general choice of the points ensures finitely many solutions of System $(*)$. The following theorem positively solves this finiteness problem.
\begin{thm}
\label{main} Under the assumption of {\rm{(i)}} and {\rm{(ii)}}, System $(*)$ consisting of the six equations \eqref{egy} has finitely many solutions.
\end{thm}
Theorem~\ref{main} improves~\cite[Theorem 1]{cifuentes20233d} where the finiteness of the solutions of System $(*)$ is only ensured for sufficiently generic choices of $ABCDEFX^*Y^*$.
In fact, let $V$ be the set of all points $A,B,C,D,E,F,X^*,Y^*$
for which either (i) or (ii) does not hold. Then the points $(A,B,C,D,E,F,X^*,Y^*)$ lie in a real affine variety $V$ in the twenty-four dimensional affine space $\mathbb{R}^{24}$; see Section~\ref{idea}. Therefore the hypothesis ``sufficiently generic choice of the eight parameters'' may be replaced by the more precise hypothesis ``choice of the corresponding points in the open set $\mathbb{R}^{24}\setminus V$''.

Another question is the unique identifiability of the solutions of System $(*)$. In this direction, we prove the following theorem.
\begin{thm}
\label{mainbis} Under the assumptions of {\rm{(i)}} and {\rm{(ii)}} on $T\in \{A,B,C,D,E,F\}$, the system consisting of the seven equations (\ref{egy}) with $T\in \{A,B,\ldots,F,G\}$ has only two solutions, namely
$(X^*,Y^*)$ and $(Y^*,X^*)$.
\end{thm}
Theorem \ref{mainbis} proves \cite[Conjecture 1]{cifuentes20233d}. The proof of Theorem \ref{main} uses computation together with some fundamental results on curves of the $3$-dimensional complex projective space. In contrast, the proof in \cite{cifuentes20233d} heavily depended on deeper results from algebraic varieties of higher dimensional complex projective space. The proof of Theorem \ref{mainbis} depends on Theorem \ref{main} and uses a projection of complex algebraic varieties.

In Section~\ref{exse}, we investigate System $(*)$ for special configurations of $ABCDEF$  satisfying (i) and (ii). We exhibit a case with infinitely many solutions and an example with more than two solutions.

The two identifiability problems above arise from inferring the 3D DNA structure of the diploid organism from partially phased data. The reader is referred to~\cite{cifuentes20233d} for details.

\section{The idea of the proof of Theorem \ref{main}}
\label{idea}
For $T\in \{A,B,\ldots,F\}$ and  $Y=(y_1,y_2,y_3)$ with $k_T\|Y-T\|^2\ne 1$, let
\begin{equation}
\label{eq241} \Gamma_T(Y)=\frac{\|Y-T\|^2}{k_T\|Y-T\|^2-1}-\|T\|^2.
\end{equation}
With this notation,  System $(*)$ can be rewritten as
\begin{equation}
\label{eq242} \|X-T\|^2=\Gamma_T(Y)+\|T\|^2, \quad T\in\{A,B,\ldots,F\}.
\end{equation}
This shows that if $(X_P,Y_P)$ is a real solution of System $(*)$ then $X_P$ is the common point of the six spheres with centers $T$ and radius $\Gamma_T(Y)+\|T\|^2$. Furthermore, if $(X_P,Y_P)$ and
$(X_Q,Y_Q)$ are two real solutions of System $(*)$ such that $X_P\neq X_Q$ but $\Gamma_T(Y_P)=\Gamma_T(Y_Q)$ for at least four $T$, say $T\in \{A,B,C,D\}$ then $A,B,C,D$ are coplanar. In fact, in this case, the four spheres of equation $\|X-T\|^2=\Gamma_T(Y_P)+\|T\|^2$  with $T\in\{A,B,C,D\}$ have two distinct common points, namely $X_P$ and $X_Q$, and hence their centers are coplanar. Therefore, an approach to the proof of Theorem \ref{main} can rely upon the study of $\Gamma_T(Y)$. We show that this approach works as we can prove that $\Gamma_T(Y)$ is constant for infinitely many $Y_P$ such that $(X_P,Y_P)$ is a solution of System $(*)$ for some $X_P$. We use assumptions (i) and (ii) in the proof. Therefore, the finiteness of the solutions of System $(*)$ holds under those two assumptions, as claimed in Theorem \ref{main}.

We also determine an affine subvariety $S$ of $\mathbb{R}^{24}$  viewed as the product of the eight $3$-dimensional real affine spaces $A,B,\ldots,F,X^*,Y^*$ such that the points $(A,B,C,D,E,F,X^*,Y^*)$ in $\mathbb{R}^{24}$ for which System $(*)$ may have infinitely many solutions lie on $S$. In $\mathbb{R}^3$, four points $A,B,C,D$ are coplanar if and only if $\Delta(A,B,C,D)=0$, see (\ref{eq04102023}). There are fifteen hypersurfaces in $\mathbb{R}^{24}$ of equations $\Delta(T_1,T_2,T_3,T_4)=0$ where $T_1,T_2,T_3,T_4$ are taken in natural alphabetic order in $\{A,B,C,D,E,F\}$.
The union of these hypersurfaces can be viewed as an affine variety $V$ of $\mathbb{R}^{24}$ whose points do not satisfy Condition (i). Similarly one can find an affine variety $W$ of $\mathbb{R}^{24}$ whose points $(A,B,C,D,E,F,X^*,Y^*)$ represent the spheres that do not satisfy Condition (ii). For this, a necessary and sufficient condition for the existence of a common (real or complex) point of four spheres is useful; see Lemma \ref{lem18set}. If $\Sigma_T$ denotes the sphere with vertex $T$ and radius $1/k_T$,  then there exists a degree six real polynomial $\Omega$ in eight indeterminates such that $\Omega(T_1,T_2,T_3,T_4,1/k_{T_1},1/k_{T_2},1/k_{T_3},1/k_{T_4})=0$ if and only if the four spheres $\Sigma_{T_i}$
have a common point. Since  $k_T$ depends on $X^*,Y^*$, substituting $k_{T_i}$ by (\ref{kxy}) gives a rational function whose numerator defines a polynomial, and hence a hypersurface in
$\mathbb{R}^{24}$. Now, $W$ can be taken to be the union of these fifteen hypersurfaces where $T_1,T_2,T_3,T_4$ is in natural alphabetic order with $T_1,T_2,T_3,T_4 \in \{A,B,C,D,E,F\}$.

Thus, $S=V\cup W$ is an affine variety of $\dim(S)<24$. Therefore, for any $(A,B,C,D,E,F,X^*,Y^*)$ from the open set $\mathbb{R}^{24}\setminus S$, System $(*)$ has finitely many solutions.

Our final remark is that assumption (i) cannot be dropped as the following case shows; see Section~\ref{exse}.

\section{Background from Euclidean geometry}
\label{back}
Our notation is standard. A point $P=(p_1,p_2,p_3)$ of $\mathbf{E}^3$ is identified with the vector $\mathbf{p}$ from to origin to $P$. Therefore, $\|P\|^2=p_1^2+p_2^2+p_3^2$, the dot product
$P\cdot Q$  equals $p_1q_q+p_2q_2+p_3q_3$, and the sphere with center $P$ and radius $r$ has equation  $\|X-P\|^2=r^2$. Furthermore, four points $A,B,C,D$ are coplanar if and only if the determinant
\begin{equation}
    \label{eq04102023}
\Delta(A,B,C,D)=\left\|
  \begin{array}{cccc}
    a_1 & b_1 & c_1 & d_1  \\
    a_2 & b_2 & c_2 & d_2  \\
    a_3 & b_3 & c_3 & d_3 \\
    1   & 1   & 1   &  1
   \end{array}
\right\|=
\left\|
  \begin{array}{cccc}
    a_1-b_1 & a_2-b_2 & a_3-b_3  \\
    a_1-c_1 & a_2-c_2 & a_3-c_3  \\
    a_1-d_1 & a_2-d_2 & a_3-d_3  \\
   \end{array}
\right\|
\end{equation}
vanishes.

We recall some classical results on configurations of spheres in the three-dimensional Euclidean space $\mathbf{E}^3$.
If $P_1,P_2$ are two distinct points of $\mathbf{E}^3$ lying on a sphere, then the center of the sphere lies on the plane passing through the middle point of the segment $\overline{P_1P_2}$ and orthogonal to the line $P_1P_2$. This indicates that if four spheres meet at least two distinct points, their centers are coplanar.

We also point out that if four spheres have at least one common complex but real point, their centers are coplanar. In fact, if $P=(p_1,p_2,p_3)$ is such a complex point, then its conjugate
$\bar{P}=(\bar{p_1}, \bar{p_2}, \bar{p_3})$ is also a common point of the same four spheres. The line $\ell=P\bar{P}$ is real; hence, it can be regarded as a line of $\mathbf{E}^3$. Also,
the middle point $M=(\ha (p_1+\bar{p_1}),\ha(p_2+\bar{p_2}),\ha(p_3+\bar{p_3}))$ lies in $\mathbf{E}^3$. Therefore, the plane passing through $M$ and orthogonal to $P\bar{P}$ is well-defined and contains the center of each of the four spheres.

\subsection{Equations linking $X$ and $Y$}\label{sec2}

\begin{lem}
\label{lemA18apr} Every solution $(X_P,Y_P)$ of System $(*)$ is uniquely determined by $Y_P$.
\end{lem}
\begin{proof}
Subtraction of the first equation from each of the other five ones gives the linear system in the coordinates of $X$
\begin{equation}
\label{eq243} 2(A-T)\cdot X=\Gamma_T(Y)-\Gamma_A(Y), \quad T\in\{B,\ldots,F\}.
\end{equation}
where
$(A-T)\cdot X=(a_1-t_1)x_1+2(a_2-t_2)x_2+2(a_3-t_3)x_3.$

Since $A,B,C,D$ are not coplanar, the determinant
$$\Delta=\left\|
  \begin{array}{cccc}
    a_1-b_1 & a_2-b_2 & a_3-b_3  \\
    a_1-c_1 & a_2-c_2 & a_3-c_3  \\
    a_1-d_1 & a_2-d_2 & a_3-d_3  \\
   \end{array}
\right\|
$$
does not vanish, and hence, Cramer's rule applies to the linear system
$$
  \begin{array}{llll}
    2(a_1-b_1)x_1+ 2(a_2-b_2)x_2 + 2(a_3-b_3)x_3 = \Gamma_B(Y_P)-\Gamma_A(Y_P);\\
    2(a_1-c_1)x_1+ 2(a_2-c_2)x_2 + 2(a_3-c_3)x_3 =\Gamma_C(Y_P)-\Gamma_A(Y_P);\\
    2(a_1-d_1)x_1 + 2(a_2-d_2)x_2 + 2(a_3-d_3)x_3 = \Gamma_D(Y_P)-\Gamma_A(Y_P),\\
  \end{array}
$$
and it gives a close formula for the coordinates of $X_P$ in functions $Y_P$ via the functions $\Gamma_T(Y)$ with $T\in\{A,B,C,D\}$.
Therefore, if $X_P=(x_1,x_2,x_3)$ then
$$x_1=\frac{1}{2\Delta} \left\|\begin{array}{cccc}
\Gamma_B(Y_P)-\Gamma_A(Y_P)  & a_2-b_2 & a_3-b_3 \\
\Gamma_C(Y_P)-\Gamma_A(Y_P) & a_2-c_2 &  a_3-c_3  \\
\Gamma_D(Y_P)-\Gamma_A(Y_P) & a_2-d_2 & a_3-d_3  \\
  \end{array}\right\|,
$$

$$x_2=\frac{1}{2\Delta} \left\|\begin{array}{ccccc}
a_1-b_1 &\Gamma_B(Y_P)-\Gamma_A(Y_P)  & a_3-b_3 \\
a_1-c_1 & \Gamma_C(Y_P)-\Gamma_A(Y_P)  & a_3-c_3  \\
a_1-d_1 &\Gamma_D(Y_P)-\Gamma_A(Y_P) & a_3-d_3  \\
  \end{array}\right\|,
$$
and
$$x_3=\frac{1}{2\Delta} \left\|\begin{array}{ccccc}
a_1-b_1 & a_2-b_2 &\Gamma_B(Y_P)-\Gamma_A(Y_P)  \\
a_1-c_1 & a_2-c_2 & \Gamma_C(Y_P)-\Gamma_A(Y_P) \\
a_1-d_1 & a_2-d_2 &\Gamma_D(Y_P)-\Gamma_A(Y_P) \\
  \end{array}\right\|
$$
which shows indeed that the coordinates of $X_P$ can uniquely be computed from those of $Y_P$ by means of $\Gamma_A(Y),\Gamma_B(Y), \Gamma_C(Y)$ and $\Gamma_D(Y)$.
\end{proof}
Now, for a solution $(X_P,Y_P)$ of System $(*)$, some equations linking the coordinates $y_1,y_2,y_3$ of $Y_P$ are established.
For this purpose, a useful fact is that for any $U=(u_1,u_2,u_3)$, the determinant
$$ \left\|\begin{array}{ccccc}
  \|U-A\|^2-\|A\|^2 & \|U-B\|^2-\|B\|^2 & \|U-C\|^2-\|B\|^2 & \|U-D\|^2-\|B\|^2 & \|U-E\|^2-\|E\|^2 \\
  a_1 & b_1 & c_1 & d_1 & e_1 \\
  a_2 & b_2 & c_2 & d_2 & e_2 \\
  a_3 & b_3 & c_3 & d_3 & e_3 \\
  1   & 1 & 1 & 1 & 1
\end{array}\right\|
$$
vanishes. Applying this to $X=(x_1,x_2,x_3)$, (\ref{eq242}) gives
\begin{equation}
\label{eqA18apr}
\left\| \begin{array}{cccccc}
  \Gamma_A(Y) & \Gamma_B(Y)& \Gamma_C(Y) & \Gamma_D(Y) & \Gamma_E(Y) \\
  a_1 & b_1 & c_1 & d_1 & e_1 \\
  a_2 & b_2 & c_2 & d_2 & e_2 \\
  a_3 & b_3 & c_3 & d_3 & e_3 \\
  1 & 1 & 1 & 1 & 1
\end{array}\right\|=0 {\mbox{ for $Y=Y_P$}}.
\end{equation}
Replacing $\{A,B,C,D,E\}$ by any five-tuple from $\{A,B,C,D,E,F\}$, disposed in natural alphabetical order, five more such determinants are obtained in the same way, say $\Phi_A(Y),\Phi_B(Y),\ldots,\Phi_F(Y)$, where the index $T$ means that $T$ is off the selected points. With this notation, the determinant in (\ref{eqA18apr}) is $\Phi_F(Y)$, and
if $(X_P,Y_P)$ is a solution of System $(*)$ then the rational functions $\Phi_T(Y)$ of $Y$  has the following property.
\begin{equation}
\label{eq16apr} {\mbox{$\Phi_T(Y_P)=0,$ for $T\in\{A,B,C,D,E,F\}$.}}
\end{equation}
Also, from the above equations,  $\|X_P-A\|^2$ can be expressed as a function of $\Gamma_A(Y_P),\Gamma_B(Y_P),\Gamma_C(Y_P)$, $\Gamma_D(Y_P)$. Thus,  (\ref{eq242}) gives
\begin{equation}
\label{eqL18apr}
\begin{array}{llll}
\vspace{0.2cm}
\Big(\left\|\begin{array}{lllll}
\Gamma_B(Y)-\Gamma_A(Y)  & a_2-b_2 & a_3-b_3 \\
\Gamma_C(Y)-\Gamma_A(Y) & a_2-c_2 &  a_3-c_3  \\
\Gamma_D(Y)-\Gamma_A(Y) & a_2-d_2 & a_3-d_3  \\
  \end{array}\right\|
  -2\Delta a_1\Big)^2+\\
  \vspace{0.2cm}
\Big(\left\|
\begin{array}{lllll}
a_1-b_1 &\Gamma_B(Y)-\Gamma_A(Y)  & a_3-b_3 \\
a_1-c_1 & \Gamma_C(Y)-\Gamma_A(Y)  & a_3-c_3  \\
a_1-d_1 &\Gamma_D(Y)-\Gamma_A(Y) & a_3-d_3  \\
\end{array}\right\|
  -2\Delta a_2\Big)^2+\\
 \vspace{0.2cm}
\Big(\left\|
\begin{array}{lllll}
a_1-b_1 & a_2-b_2 &\Gamma_B(Y)-\Gamma_A(Y)  \\
a_1-c_1 & a_2-c_2 & \Gamma_C(Y)-\Gamma_A(Y) \\
a_1-d_1 & a_2-d_2 &\Gamma_D(Y)-\Gamma_A(Y) \\
\end{array}\right\|
  -2\Delta a_3\Big)^2-\\
\vspace{0.2cm}
  4\Delta^2(\Gamma_A(Y)+\|A\|^2)=0,\,\,{\mbox{for $Y=Y_P$}}.
\end{array}
\end{equation}
Replacing $\{A,B,C,D\}$ by any four-tuple from $\{A,B,C,D,E,F\}$, disposed in natural alphabetical order, fourteen more such functions of $Y$ arise in the same way.  These fifteen functions are denoted by  $\Psi_{A,B}(Y),\Psi_{A,C}(Y),\ldots,$ $\Psi_{E,F}(Y)$, where the index ${T_1,T_2}$ means that $T_1$ and $T_2$ are off the selected points. With this notation, the function in (\ref{eqL18apr}) is $\Psi_{E,F}(Y)$. Therefore, if $Y_P$ is a solution of System $(*)$ then
\begin{equation}
\label{eqA17apr} {\mbox{$\Psi_{T_1,T_2}(Y_P)=0$ for any $T_1,T_2\in\{A,B,C,D,E,F\}$ }}
\end{equation}
where $T_1$ precedes $T_2$ in the natural alphabetic order.
\subsection{Four spheres  with a common point}
The arguments in subsection \ref{sec2} are also useful to prove the following lemma.
\begin{lem}
\label{lem18set} For $i=1,2,3,4$, let $\Sigma_i$ be the sphere of center $P_i=(p_{i1}, p_{i2},p_{i3})$ and radius $r_i$. Then $\Sigma_1,\Sigma_2,\Sigma_3,\Sigma_4$ have a common real or complex point if and only if $\Omega(P_1,P_2,P_3,P_4,r_1^2,r_2^2,r_3^2,r_4^2)=$
$$
\begin{array}{llll}
\vspace{0.2cm}
\Big(\left\|\begin{array}{lllll}
(r_2^2-r_1^2)-(\|P_2\|^2-\|P_1\|^2)  & p_{12}-p_{22} & p_{13}-p_{23} \\
(r_3^2-r_1^2)-(\|P_3\|^2-\|P_1\|^2)  & p_{12}-p_{32} & p_{13}-p_{33}  \\
(r_4^2-r_1^2)-(\|P_4\|^2-\|P_1\|^2)  & p_{12}-p_{42} & p_{13}-p_{43}  \\
  \end{array}\right\|
  -2
  \left\|
  \begin{array}{cccc}
    p_{11}-p_{21} & p_{12}-p_{22} & p_{13}-p_{23}  \\
    p_{11}-p_{31} & p_{12}-p_{32} & p_{13}-p_{33}  \\
    p_{11}-p_{41} & p_{12}-p_{42} & p_{13}-p_{43}  \\
   \end{array}
\right\|
  p_{11}\Big)^2+\\
  \vspace{0.2cm}
\Big(\left\|
\begin{array}{lllll}
p_{11}-p_{21} & (r_2^2-r_1^2)-(\|P_2\|^2-\|P_1\|^2) & p_{13}-p_{23} \\
p_{11}-p_{31} & (r_3^2-r_1^2)-(\|P_3\|^2-\|P_1\|^2) & p_{13}-p_{33}  \\
p_{11}-p_{41} & (r_4^2-r_1^2)-(\|P_4\|^2-\|P_1\|^2) & p_{13}-p_{43}  \\
\end{array}\right\|
  -2\left\|
  \begin{array}{cccc}
    p_{11}-p_{21} & p_{12}-p_{22} & p_{13}-p_{23}  \\
    p_{11}-p_{31} & p_{12}-p_{32} & p_{13}-p_{33}  \\
    p_{11}-p_{41} & p_{12}-p_{42} & p_{13}-p_{43}  \\
  \end{array}
\right\|p_{12}\Big)^2+\\
 \vspace{0.2cm}
\Big(\left\|
\begin{array}{lllll}
p_{11}-p_{21} & p_{12}-p_{22} &  (r_2^2-r_1^2)-(\|P_2\|^2-\|P_1\|^2) \\
p_{11}-p_{31} & p_{22}-p_{32} & (r_3^2-r_1^2)-(\|P_3\|^2-\|P_1\|^2) \\
p{11}-p_{41} & p_{22}-p_{42} &  (r_4^2-r_1^2)-(\|P_4\|^2-\|P_1\|^2) \\
\end{array}\right\|
  -2\left\|
  \begin{array}{cccc}
    p_{11}-p_{21} & p_{12}-p_{22} & p_{13}-p_{23}  \\
    p_{11}-p_{31} & p_{12}-p_{32} & p_{13}-p_{33}  \\
    p_{11}-p_{41} & p_{12}-p_{42} & p_{13}-p_{43}  \\
   \end{array}
\right\| p_{13}\Big)^2-\\
\vspace{0.2cm}
  4\Big(\left\|
  \begin{array}{cccc}
  p_{11}-p_{21} & p_{12}-p_{22} & p_{13}-p_{23}  \\
    p_{11}-p_{31} & p_{12}-p_{32} & p_{13}-p_{33}  \\
    p_{11}-p_{41} & p_{12}-p_{42} & p_{13}-p_{43}  \\
  \end{array}
\right\|^2+\|P_1\|^2+r_1^2\Big)=0.
\end{array}
$$
\end{lem}
\begin{proof} The four spheres $\Sigma_i$ have a common point if and only if the system consisting of the four equations $\|X-P_i\|^2=r_i^2$ has a solution. For $i=1,2,3,4$, let $\Lambda_i=r_i^2-\|P_i\|^2$. Then this system becomes $ \|X\|^2-2X\cdot P_i -\Lambda_i=0,\,\, i\in \{1,2,3,4\}.$ Subtraction of the first equation from the other three gives
$$2 (P_1-P_i)\cdot X=\Lambda_i-\Lambda_1,\,\, i\in \{2,3,4\}.$$
Comparison with (\ref{eq243}) shows that the arguments in subsection \ref{sec2} can be used for the present proof. In particular, (\ref{eqL18apr}) holds if $(A,B,C,D)$ is replaced by $(P_1,P_2,P_3,P_4)$, and $(\Gamma_A(Y),\Gamma_B(Y),\Gamma_C(Y),\Gamma_D(Y))$ by $(\Lambda_1,\Lambda_2,\Lambda_3,\Lambda_4)$. Since $\Lambda_i=r_i^2-\|P_i\|^2$, the claim follows.
\end{proof}

\subsection{Basic facts on surfaces and algebraic curves in the $3$-dimensional complex space.}
\label{bfa}
In the $3$-dimensional affine space $\mathbb{C}^3$ over the complex field $\mathbb{C}$, an affine variety (or affine set) $W$ consists of all common points of a finite number of surfaces. If $f_1(X,Y,Z),\ldots, f_k(X,Y,Z)$ are the polynomial associated to those surfaces then the points of $V$ viewed as ordered triples $(p_1,p_2,p_3)$ are exactly the solutions of the system $f_1(X,Y,Z)=0,\ldots, f_k(X,Y,Z)=0$. An affine variety $W$ is \emph{reducible} if there exist affine varieties $W_1\subset W$, $W_2\subset W$, $W_1\neq W_2$ such that $W=W_1\cup W_2$. Otherwise, $W$ is \emph{irreducible}. Every affine variety is the union of finitely many irreducible ones. Since we are in dimension three, $W$ is the union of a finite number of affine irreducible varieties of dimensions $0$, $1$, and $2$. Assume that $W$ has infinitely many real points, and take any infinite subset $G$ of real points of $W$. Since any zero-dimensional affine variety consists of a finite number of points, $W$ contains either an irreducible curve or an irreducible surface $V$, which still contains infinitely many real points contained in $G$. Take any finite subset $F$ of $G$ in the latter case. From a classical theorem going back to Bertini, see \cite[Theorem 1.8]{charles2016bertini}, there exists an irreducible curve lying on $V$ that passes through every point in $F$. In \cite[Theorem 1.8]{charles2016bertini} states: Let $\cX$ be an irreducible variety of dimension $m \ge 2$ over a field $k$. Let $F \subset  \cX$ be a finite set of closed points. Then, a geometrically irreducible variety exists $\cY \subseteq \cX$ of dimension $m-1$ containing $F$. We will apply this result to the case where $m=2$, $k=\mathbb{C}$, and $|F|=1$. Since $\mathbb{C}$ is algebraically closed, the closed points of $\cX$ are exactly the points of $\cX$. For $k=\mathbb{R}$, a closed point is either a real point or a pair of two points, say $P,Q$, such that $Q$ is the conjugate of $P$ (and $P$ is the conjugate of $Q$).

For the study of an irreducible projective complex algebraic curve $\cC$ of the $3$-dimensional complex projective space, we use concepts and results about its branches (places), as well as about the field $\mathbb{C}((\tau))$ of formal power series $\alpha(\tau)=\sum_{i\ge \ell}^{\infty}u_i\tau^i$ in an indeterminate $\tau$. Here, $\mathbb{C}((\tau))$ is the rational field of the ring $\mathbb{C}[[\tau]]$ of all formal power series $\alpha(\tau)=\sum_{i\ge 0}^{\infty}u_i\tau^i$.
Our notation and terminology are standard; see \cite{lefschetz2012algebraic,shafaverich2013basic,seidenberg1968elements}. In particular, $\alpha=0$ if and only if $u_i=0$ for every $i$, and for $\alpha\neq 0$, ${\rm{ord}}(\alpha(\tau))$ denotes the smallest $i$ with $u_i\neq 0$.

Every point $P\in \cC$ is the center of at least one branch of $\cC$. If $P$ is a non-singular point of $\cC$, then exactly one branch of $\cC$ is centered at $P$. For a singular point of $\cC$ there may be more then one but finitely many branches centered at that point. If $P=(\xi,\eta,\zeta)=(\xi:\eta:\zeta:1)\in\cC$ is an affine point, a primitive representation of a branch $\gamma$ of $\cC$ centered at $P$ is a triple $(x(\tau),y(\tau),z(\tau))$ so that 
\begin{equation}
\label{eqA20apr}
x(\tau)=\xi+\varphi_1(\tau),\, y(\tau)=\eta+\varphi_2(\tau),\,z(\tau)=\zeta+\varphi_3(\tau)
\end{equation}
with $\varphi_i(\tau)\in \mathbb{C}[[\tau]]$ where at least one $\varphi_i(\tau)$ does not vanish.

Let $P=(\eta_1:\eta_2:\eta_3:0)$ be a point of $\cC$ at infinity. A primitive branch representation gives a branch of $\cC$ centered at $P$
$(x(\tau):y(\tau):z(\tau):1)$ with $x(\tau),y(\tau),z(\tau)\in \mathbb{C}((\tau))$ of the form
$$x(\tau)=(\eta_1+\varphi_1(\tau))/\alpha(\tau), y(\tau)=(\eta_2+\varphi_2(\tau))/\alpha(\tau),z(\tau)=(\eta_3+\varphi_3(\tau))/\alpha(\tau)$$  where $\varphi_i(\tau)\in \mathbb{C}[[\tau]]$, $\alpha(\tau)=t^m \beta(\tau)$ with $m\ge 1$ and $\min_i\{{\rm{ord}}(\varphi_i(\tau))\}<m$.

A branch representation is primitive (or irreducible) if the ${\rm{g.c.d.}}$ of all exponents of $\tau$ in the power series $\varphi_i(\tau)$ is equal to $1$.

Let $M(Y)$ be a polynomial in $Y=(y_1,y_2,y_3)$. Then $M(Y)$ is said to vanish at a branch of $\cC$ given by a primitive representation $(x(\tau),y(\tau),z(\tau))$ if $M(x(\tau),y(\tau),z(\tau))$ is the zero power series. If $M(Y)$ vanishes at some branch of $\cC$, it vanishes at every branch of $\cC$ as $\cC$ is irreducible.
Let $L(Y)=L_1(Y)/L_2(Y)$ be a complex rational function with polynomials $L_1(Y),L_2(Y)$ which do not vanishes at any branch of $\cC$.
For a branch $\gamma$  of $\cC$ given by a primitive representation $(x(\tau),y(\tau),z(\tau))$, let $L_j(\tau)=L_j(x(\tau),y(\tau),z(\tau))$ with $j=1,2$. Then $L_j(\tau)\in \mathbb{C}((\tau))$, and $\gamma$ is a pole of $L(Y)$ if ${\rm{ord}}(L_1(\tau))<{\rm{ord}}(L_2(\tau))$, and it is a zero of $L(Y)$ if
${\rm{ord}}(L_1(\tau))> {\rm{ord}}(L_2(\tau))$.
By a classical result, see \cite[Theorem 14.1]{seidenberg1968elements}, every complex rational function $L$
has only finitely many poles, and $L$ has some
poles, unless $L$ is a (non-zero) constant.
\subsection{Proof of Theorem \ref{main}}
From Lemma \ref{lemA18apr}, the number of solutions $(X_P,Y_P)$ of System $(*)$ can be obtained by counting their component $Y_P$. If this number is infinite as in the hypothesis of Theorem \ref{main},
there are infinitely many such $Y_P$'s satisfying both (\ref{eq16apr}) and (\ref{eqA17apr}). This suggests performing a proof of Theorem \ref{main} based upon the concepts and results from algebraic geometry recalled in subsection \ref{bfa}.

For this purpose, if the system involves rational functions as in (\ref{eq16apr}) and (\ref{eqA17apr}), the usual approach is to take the numerators of the rational functions and investigate the
(possible reducible) affine variety defined by those numerators.
Let $\cN(\Phi_T(Y))$  and $\cN(\Psi_{T_1,T_2}(Y))$ denote, respectively,
the numerators of the rational functions $\Phi_T(Y)$ and $\Psi_{T_1,T_2}(Y)$ introduced in Section \ref{sec2}.
In particular,
\begin{equation}
\label{eq16set}
\cN(\Phi_F(Y))=\left\| \begin{array}{cccccc}
  \Gamma_A(Y) & \Gamma_B(Y)& \Gamma_C(Y) & \Gamma_D(Y) & \Gamma_E(Y) \\
  a_1 & b_1 & c_1 & d_1 & e_1 \\
  a_2 & b_2 & c_2 & d_2 & e_2 \\
  a_3 & b_3 & c_3 & d_3 & e_3 \\
  1 & 1 & 1 & 1 & 1
\end{array}\right\|
\prod_{T\in\{A,B,C,D,E\}}(k_T\|Y-T\|^2-1),
\end{equation}
and
$$
\cN(\Psi_{E,F}(Y))=
$$
\begin{equation}
\label{eqL18aprA}
\begin{array}{llll}
\vspace{0.2cm}
\Big[\Big(\left\|\begin{array}{lllll}
\Gamma_B(Y)-\Gamma_A(Y)  & a_2-b_2 & a_3-b_3 \\
\Gamma_C(Y)-\Gamma_A(Y) & a_2-c_2 &  a_3-c_3  \\
\Gamma_D(Y)-\Gamma_A(Y) & a_2-d_2 & a_3-d_3  \\
  \end{array}\right\|
  -2\Delta a_1\Big)\Big(\prod_{T\in\{A,B,C,D\}}(k_T\|Y-T\|^2-1)\Big)\Big]^2+\\
  \vspace{0.2cm}
\Big[\Big(\left\|
\begin{array}{lllll}
a_1-b_1 &\Gamma_B(Y)-\Gamma_A(Y)  & a_3-b_3 \\
a_1-c_1 & \Gamma_C(Y)-\Gamma_A(Y)  & a_3-c_3  \\
a_1-d_1 &\Gamma_D(Y)-\Gamma_A(Y) & a_3-d_3  \\
\end{array}\right\|
  -2\Delta a_2\Big)\Big(\prod_{T\in\{A,B,C,D\}}(k_T\|Y-T\|^2-1)\Big)\Big]^2+\\
 \vspace{0.2cm}
\Big[\Big(\left\|
\begin{array}{lllll}
a_1-b_1 & a_2-b_2 &\Gamma_B(Y)-\Gamma_A(Y)  \\
a_1-c_1 & a_2-c_2 & \Gamma_C(Y)-\Gamma_A(Y) \\
a_1-d_1 & a_2-d_2 &\Gamma_D(Y)-\Gamma_A(Y) \\
\end{array}\right\|
  -2\Delta a_3\Big)\Big(\prod_{T\in\{A,B,C,D\}}(k_T\|Y-T\|^2-1)\Big)\Big]^2-\\
\vspace{0.2cm}
  \Big[4\Delta^2(\Gamma_A(Y)+\|A\|^2)\Big]\Big[\prod_{T\in\{A,B,C,D\}}(k_T\|Y-T\|^2-1)\Big]^2.
\end{array}
\end{equation}

Let $W$ be the complex affine variety associated with $\cN(\Phi_T(Y))$ and $\cN(\Psi_{T_1,T_2}(Y))$
where $T,T_1,T_2\in \{A,B,C,D,E,F\}$ and $T_1$ precedes $T_2$ in the natural alphabetic order.
Obviously, $\dim(W)\le 2$. As pointed out in subsection~\ref{bfa}, we may assume that $W$ contains an (affine) irreducible curve $\cD$ through a point $Y_P$ such that $(X_P,Y_P)$ is a solution of System $(*)$. In particular, $Y_P$ satisfies   both (\ref{eq16apr}) and (\ref{eqA17apr}). Since $\cD$ is contained in $W$, for every point $Y\in\cD$ we have $\Phi_T(Y)=0$ and $\Psi_{T_1,T_2}(Y)=0$ for $T,T_1,T_2\in\{A,B,C,D,E,F\}$ where $T_1$ precedes $T_2$ in the natural alphabetic order.

Let $\cC$ be the projective closure of $\cD$. Then $\cC$  is a complex, projective, irreducible algebraic curve of $P_3(\mathbb{C})$ whose affine points lie on $W$.

\begin{lem}
\label{lem20apr} For any $T\in\{A,B,C,D,E,F\}$, neither $\|Y-T\|^2$ nor $k_T\|Y-T\|^2-1$ vanishes at every branch of $\cC$.
\end{lem}
\begin{proof}  As we have already observed, there exists $Y_P=(\eta_1:\eta_2:\eta_3:1)\in \cC$ such that $(X_P,Y_P)$ is a solution of System $(*)$. Let $\gamma$ be a branch of $\cC$ centered at $Y_P$. Since $Y_P$ is an affine point, $\gamma$ has a primitive branch representation (\ref{eqA20apr}).
Assume on the contrary $k_T\|Y-T\|^2=1$ at every branch of $\cC$. Then $k_T((x(\tau)-t_1)^2+(y(\tau)-t_2)^2+(z(\tau)-t_3)^2)=1$ whence
$k_T((\eta_1-t_1)^2+(\eta_2-t_2)^2+(\eta_3-t_3)^2+\psi(\tau))=1$ where $\psi\in \mathbb{C}[[\tau]]$ with ${\rm{ord}}(\psi(\tau))>0$. Therefore, $k_T((\eta_1-t_1)^2+(\eta_2-t_2)^2+(\eta_3-t_3)^2)=1$,
that is, $k_T(\|Y_P-T\|^2)-1=0$. But then there exists no $X_P$ such that $(X_P,Y_P)$ is a solution of $(*)$, a contradiction. The proof for the case $\|Y-T\|^2=0$ is analogous.
\end{proof}
By Lemma \ref{lem20apr}, $\Gamma_T(Y)$ does not vanish at any branch of $\cC$ for $T\in\{A,B,\ldots,F\}$.
\begin{lem}
\label{lem17apr} If the branch $\gamma$ of $\cC$ is centered at an affine point, then $\gamma$ is not a pole of $\Gamma_T(Y)$ with $T\in \{A,B,C,D,E,F\}$.
\end{lem}
\begin{proof} By way of a contradiction, assume that a branch $\gamma$ is a pole of $\Gamma_A(Y)$. Let $\eta_P$ be the center.  From (\ref{eq16set}),
 the coefficient of $\Gamma_T(Y)$ in $\Phi_F(Y)$ with $T\neq F$ does not vanish as no four points from $A,B,C,D,E$ are coplanar. Since ${\rm{ord}}(\Gamma_A(\gamma))<0$ and $\Phi_F(\gamma)=0$, this yields that ${\rm{ord}}(\Gamma_T(\gamma))<0$  for another $T\in\{B,C,D,E\}$. W.l.g., we may suppose this occurs for $\Gamma_B(Y)$. Then $\gamma$ is also a pole of $\Gamma_B(Y)$. Repeat the same argument for $\Gamma_B(Y)$ and $\Phi_A(Y)$. It turns out that $\gamma$ is also a pole of one of the function $\Gamma_T(Y)$ with $T\in \{C,D,E,F\}$, say $T=C$. We show that the $\gamma$ is not a pole of any of the remaining functions $\Gamma_T(Y)$ with $T\in \{D,E,F\}$.
We may  assume on the contrary that $\gamma$ is a pole of each $\Gamma_T(Y)$ with $T\in\{A,B,C,D\}$.
Therefore $\gamma$ is also a pole of $\Gamma_T(Y)-\|T\|^2$ with $T\in \{A,B,C,D\}$.
Then $\gamma$ is a zero of $(\Gamma_T(Y)-\|T\|^2)^{-1}$, and hence $\|\eta_P-T\|^2-1/k_T=0$ for\ $T\in \{A,B,C,D\}$.
This means that $\eta_P$ is a common point of four real spheres, namely those of center $T$ and radius $1/\sqrt{k_T}$ where $T\in \{A,B,C,D\}$.
Observe that $\eta_P$ is not a complex point by the final remark in subsection \ref{back}. Thus, $\eta_P$ is a real point, which contradicts assumption (ii).

Now, replace $(A,B,C)$ by $(D,E,F)$. Then $\gamma$ is a pole for $\Gamma_D(Y)$ but it is not for any $\Gamma_T(Y)$ with $T\in\{A,B,C\}$, that is,
$\Gamma_D(\gamma)=ct^{-i}+\ldots$ with $c\neq 0, i>0$, whereas ${\rm{ord}}(\Gamma_T(\gamma))\ge 0$ for $T\in\{A,B,C\}$. Furthermore, $\Phi_{E,F}(\gamma)=0$.
In particular, the coefficient of $t^{-2i}$ in $\Gamma_{E,F}(\gamma)$ must vanish.

On the other hand, (\ref{eqL18apr}) shows that in $\Phi_{E,F}$ the degree of  $\Gamma_D(Y)$ is two,
and the coefficient of $t^{-2i}$ in $\Gamma_D(Y)^2(\gamma)$  is
\begin{equation}
\label{eqT19apr}
\left\| \begin{array}{lllll}
a_2-b_2 & a_3-b_3  \\
a_2-c_2 & a_3-c_3
\end{array}
\right \|^2+
\left\| \begin{array}{lllll}
a_1-b_1 & a_3-b_3  \\
a_1-c_1 & a_3-c_3
\end{array}
\right \|^2+
\left\| \begin{array}{lllll}
a_1-b_1 & a_2-b_2  \\
a_1-c_1 & a_2-c_2
\end{array}
\right \|^2.
\end{equation}
Since this coefficient must vanish,
\begin{equation}
\label{eqT19aprbis}
\left\| \begin{array}{lllll}
a_2-b_2 & a_3-b_3  \\
a_2-c_2 & a_3-c_3
\end{array}
\right \|=
\left\| \begin{array}{lllll}
a_1-b_1 & a_3-b_3  \\
a_1-c_1 & a_3-c_3
\end{array}
\right \|=
\left\| \begin{array}{lllll}
a_1-b_1 & a_2-b_2  \\
a_1-c_1 & a_2-c_2
\end{array}
\right \|=0.
\end{equation}
This implies that $A,B, and C$ are collinear, contradicting assumption (i).
Since $\Psi_{E,F}(Y)=0$,  we have that $\gamma$ must also  be a pole of one of the functions $\Gamma_A(Y),\Gamma_B(Y)$, $\Gamma_C(Y)$, a contradiction.
\end{proof}
 Next, we show that if $P=(\eta_0:\eta_1:\eta_2:0)$ is a point of $W$ at infinity, then $\eta_1^2+\eta_2^2+\eta_3^2=0$. For this, we need to use
 homogeneous coordinates $(y_1,y_2,y_3,t)$. With $$\tilde{\Gamma}_T(Y)=\frac{\|Y-tT\|^2}{k_T\|Y-tT\|^2-t^2}-\|tT\|^2,$$  Equations (\ref{eq16set}) becomes
\begin{equation}
\label{eq16Bset}
\cN(\Phi_F(Y))=\left\| \begin{array}{cccccc}
\tilde{\Gamma}_A(Y) & \tilde{\Gamma}_B(Y)& \tilde{\Gamma}_C(Y) & \tilde{\Gamma}_D(Y) & \tilde{\Gamma}_E(Y) \\
  a_1 & b_1 & c_1 & d_1 & e_1 \\
  a_2 & b_2 & c_2 & d_2 & e_2 \\
  a_3 & b_3 & c_3 & d_3 & e_3 \\
  1 & 1 & 1 & 1 & 1
\end{array}\right\|
\prod_{T\in\{A,B,C,D,E\}}(k_T\|Y-tT\|^2-t^2)
\end{equation}

This shows that the point $P=(\eta_0:\eta_1:\eta_2:0)$ lies on the projective closure $\mathcal{F}$ of the hypersurface associated with $\cN(\Phi_F(Y))$ if and only if
\begin{equation}
\label{eq16setC}
(\eta_1^2+\eta_2^2+\eta_3^2)^4
\left\| \begin{array}{cccccc}
k_Bk_Ck_Dk_E &  k_Ak_Ck_Dk_E & k_Ak_Bk_Dk_E & k_Ak_Bk_Ck_E & k_Ak_Bk_Ck_D  \\
  a_1 & b_1 & c_1 & d_1 & e_1 \\
  a_2 & b_2 & c_2 & d_2 & e_2 \\
  a_3 & b_3 & c_3 & d_3 & e_3 \\
  1 & 1 & 1 & 1 & 1
\end{array}\right\|=0
\end{equation}
Therefore, if $\eta_1^2+\eta_2^2+\eta_3^2\neq 0$, then the determinant in (\ref{eq16setC}) vanishes, and hence every point at infinity is a point of $\mathcal{F}$. However, this is a contradiction as the projective closure of the affine hypersurface cannot contain all points at infinity.

\begin{lem}
\label{lem16apr} If the branch $\gamma$ is centered at a point of\, $\cC$ at infinity, then $\gamma$ is not a pole of $\Gamma_T(Y)$ for any $T\in \{A,B,C,D,E,F\}$.
\end{lem}
\begin{proof} Let $P=(\eta_1:\eta_2:\eta_3:0)$ be a point of $\cC$ at infinity.
 As we have already observed, $\eta_1^2+\eta_2^2+\eta_3^2=0$. Let $\gamma$ be a branch of $\cC$ centered at $P$ with a primitive branch representation
$(x(\tau):y(\tau):z(\tau):1)$ where
$$x(\tau)=(\eta_1+\varphi_1(\tau))/\alpha(\tau), y(\tau)=(\eta_2+\varphi_2(\tau))/\alpha(\tau),z(\tau)=(\eta_3+\varphi_3(\tau))/\alpha(\tau)$$  where $\alpha(\tau)=\tau^{-m} \beta(\tau)$ with $m\ge 1$ and $\beta(\tau)\in \mathbb{C}[[\tau]],  \beta(0)\neq 0$.

Suppose that  $\gamma$ is a pole of $\Gamma_A(Y)$.  Then $\gamma$ is also a pole of $\Gamma_A(Y)+\|A\|^2$. Hence $\gamma$ is a zero of $(\Gamma_A(Y)+\|A\|^2)^{-1}$ and, that is, $\gamma$ is a zero of
\begin{equation}
\label{eqF17apr}
k_A-\frac{1}{\|Y+A\|^2}.
\end{equation}
 A straightforward computation shows
$$
\|Y+A\|^2=(x(\tau)+a_1)^2+(y(\tau)+a_2)^2+(z(\tau)-a_3)^2=\\
$$
$$
\frac{\sum_{i=1}^3(\eta_i+\varphi_i(\tau)+a_i)^2}{\alpha(\tau)^2}=\frac{\sum_{i=1}^3 \eta_i^2+2\sum_{i=1}^3 \eta_ia_i+2\sum a_i\varphi_i(\tau)+\|A\|^2+2\sum_{i=1}^3 \eta_i\varphi_i(\tau)
+\sum_{i=1}^3 \varphi_i^2(\tau)}{\alpha(\tau)^2}.
$$

Expanding (\ref{eqF17apr}) gives
$$k_A- \frac{ 1}{U_1(\tau)+U_2(\tau)+U_3(\tau)+V_1(\tau)+V_2(\tau)}$$
with
$$
\begin{array}{llll}
U_1(\tau)=2(\eta_1a_1+\eta_2a_2+\eta_3a_3)\tau^{-m}\beta(\tau)^{-1},\\
U_2(\tau)=2(a_1\varphi_1(\tau)+a_2\varphi_2(\tau)+a_3\varphi_3(\tau))\tau^{-m}\beta(\tau)^{-1},\\
U_3(\tau)=\|A\|^2,
\end{array}
$$
and
$$
\begin{array}{lll}
V_1(\tau)=2(\eta_1\varphi_1(\tau)+\eta_2\varphi_2(\tau)+\eta_3\varphi_3(\tau))\tau^{-2m}\beta(\tau)^{-2},\\
V_2(\tau)=(\varphi_1(\tau)^2+\varphi_2(\tau)^2+\varphi_3(\tau)^2)\tau^{-2m}\beta(\tau)^{-2}.
\end{array}
$$
If $a_1\eta_1+a_2\eta_2+a_3\eta_3\neq 0$ then
$${\rm{ord}}(U_1(\tau))=-m<1-m\le{\rm{ord}}(2(a_1\varphi_1(\tau)+a_2\varphi_2(\tau)+a_3\varphi_3(\tau))-m\leq {\rm{ord}}(U_2(\tau)),$$ and ${\rm{ord}}(U_1(\tau))=-m<0={\rm{ord}}(U_3(\tau))$ whence ${\rm{ord}}(U_1(\tau))={\rm{ord}}(U_1(\tau))+U_2(\tau)+U_3(\tau)).$  Since $\beta(0)\neq 0$ and $\gamma$ is not a pole of $(\Gamma_A(Y)+\|A\|^2)^{-1}$
otherwise $\gamma$ would be a zero of $\Gamma_A(Y)+\|A\|^2$ and hence $\gamma$ would not be a pole of $\Gamma_A(Y)$),
This implies that ${\rm{ord}}(U_1(\tau))={\rm{ord}}(V_1(\tau)+V_2(\tau))$. Hence
$2(a_1\eta_1+a_2\eta_2+a_3\eta_3)\beta(0)^{-1}$ equals the opposite of the coefficient of $\tau^{-m}$ in $V_1(\tau)+V_2(\tau)$. Since the latter power series does not depend on $A$, there is a constant $\kappa$ only depending on $\gamma$ such that
\begin{equation}
\label{eqE17apr} a_1\eta_1+a_2\eta_2+a_3\eta_3=\kappa.
\end{equation}
As in the proof of Lemma \ref{lem17apr}, it is possible to prove that $\gamma$ is also a pole of two more functions $\Gamma_T(Y)$, say $\Gamma_B(Y)$ and $\Gamma_C(Y)$. Since (\ref{eqE17apr}) holds true if $A$ is replaced either by $B$ or by $C$, we obtain
$$
\begin{array}{llll}
a_1\eta_1+a_2\eta_2+a_3\eta_3=\kappa,\\
b_1\eta_1+b_2\eta_2+b_3\eta_3=\kappa,\\
c_1\eta_1+c_2\eta_2+c_3\eta_3=\kappa,\\
\eta_1^2+\eta_2^2+\eta_3^2=0.
\end{array}
$$
From the first three equations,
$$\eta_1=\kappa \left\|
  \begin{array}{cccc}
    1 &  a_2& a_3   \\
    1 &  b_2& b_3   \\
    1 & c_2& c_3   \\
    \end{array}
\right\|\Delta^{-1},
\eta_2=\kappa \left\|
  \begin{array}{cccc}
    a_1 & 1 & a_3   \\
    b_1 & 1 & b_3   \\
    c_1 & 1& c_3   \\
    \end{array}
\right\|\Delta^{-1},
\eta_3=\kappa \left\|
  \begin{array}{cccc}
    a_1 &  a_2& 1   \\
    b_1 &  b_2& 1   \\
    c_1 & c_2& 1   \\
    \end{array}
\right\|\Delta^{-1}
$$
where
$$\Delta=\left\|
  \begin{array}{cccc}
    a_1 & a_2 & a_3   \\
    b_1 & b_2 & b_3   \\
    c_1 & c_2 & c_3   \\
    \end{array}
\right\|.
$$
Therefore,
$$\eta_1^2+\eta_2^2+\eta_3^2= \kappa^2 \left(
 \left\|
  \begin{array}{cccc}
    a_1 & 1 & a_3   \\
    b_1 & 1 & b_3   \\
    c_1 & 1& c_3   \\
    \end{array}
\right\|^2+\left\|
  \begin{array}{cccc}
    1 & a_2 & a_3   \\
    1 & b_2 & b_3   \\
    1 & c_2 & c_3   \\
    \end{array}
\right\|^2+\left\|
  \begin{array}{cccc}
    a_1 &  a_2& 1   \\
    b_1 &  b_2& 1   \\
    c_1 & c_2& 1   \\
    \end{array}
\right\|^2
\right)\Delta^{-2}.
$$

By the fourth equation, this implies $\eta_1=\eta_2=\eta_3=0$. Therefore (\ref{eqT19aprbis}) holds and hence $A,B,C$ are collinear points, a contradiction.
\end{proof}
Lemmas \ref{lem17apr} and \ref{lem16apr} have the following corollary.
\begin{lem}
\label{prop19apr} $\Gamma_T(Y)$ with $T\in \{A,B,C,D,E,F\}$ is a constant on the branches of $\cC$.
\end{lem}
From Lemma \ref{prop19apr}, there exists $c_T\in \mathbb{C}$ such that $\Gamma_T(Y)-\|T\|^2=c_T$ for all  affine points $Y\in \cC$.
Take two of them, say $Y_P$ and $Y_Q$ such that both $(X_P,Y_P)$ and $(X_Q,Y_P)$ are real solutions of System $(*)$ for some $X_P$ and $X_Q$. Then (\ref{eq242}) yields $\|X_P-T\|^2=c_T=\|X_Q-T\|^2$ for $T\in \{A,B,C,D,E,F\}$.
Therefore, $X_P$ and $X_Q$ are common points of six spheres with centers at the points $A,B,\ldots,F$ contradicting assumption (i).
This completes the proof of Theorem~\ref{main}.

\subsection{Special configurations}
\label{exse}
For $A,B,C,D,E,F$ choose any six points lying on a circle $\Sigma$ of $\mathbb{E}^3$, and for the $X^*,Y^*$ take two points on the line $\ell$  through the center of $\Sigma$ and
orthogonal to the plane containing $\Sigma$. For any two points $X,Y\in \ell$, if $T$ ranges over $\{A,B,C,D,E,F\}$, then $1/\|X-T\|^2  + 1/\|Y-T\|^2$ does not change. Now, take any
$X\in\ell$ with $$\frac{1}{\|X-A\|^2} <\frac{1}{\|X^*-A\|^2} +\frac{1}{\|Y^*-A\|^2}.$$ Then there exists $Y\in \ell$  such that $(X,Y)$ is a solution of System $(*)$. We have infinitely many such pairs $(X,Y)$.
However,  some tests made by {\bf{Macaualy2}} show that if $X\not\in \ell$, then System $(*)$ may only have finitely many solutions. It remains open the problem whether Condition (i) alone is enough to ensure the finiteness of the solutions of System $(*)$.

A generalization of the above construction shows that Conditions (i) and (ii) are not sufficient to ensure unique identifiability for System $(*)$. This time, take two circles $\Sigma_1$ and $\Sigma_2$ in $\mathbf{E}^3$ lying on two parallel planes $\pi_1$ and $\pi_2$ such that the line $\ell$ joining their centers be orthogonal to the planes. For $A,B,C$, take any three points on $\Sigma_1$, and for $D,E,F$, any three points on $\Sigma_2$. Furthermore, let $X^*,Y^*$ be any two distinct points on $\ell$.  Fix a reference system $(x,y,z)$ in $\mathbf{E}^3$ so that $\ell$ is the $z$-axis. Then $\pi_1$ and $\pi_2$ have equations $z=a_3$ and $z=d_3$, respectively. Furthermore, $A=(a_1,a_2,a_3),B=(b_1,b_2,a_3),C=(c_1,c_2,a_3)$ with $a_1^2+a_2^2=b_1^2+b_2^2=c_1^2+c_2^2=\rho_1^2$. Let $\rho_1^2+a_3^2=r_1^2$. Then $\overline{OA}^2=\overline{OB}^2=\overline{OC}^2=r_1^2$. Similarly, for $D=(d_1,d_2,d_3),E=(e_1,e_2,d_3),F=(f_1,f_2,d_3)$, $d_1^2+d_2^2=e_1^2+e_2^2=f_1^2+f_2^2=\rho_2^2$. Let $\rho_2^2+d_3^2=r_2^2$. Moreover, the points on $\ell$ are parametrized by a unique variable. Hence, we may assume $X=(0,0,U)$ and $Y=(0,0,V)$, in particular we set $X^*=(0,0,u)$, $Y^*(0,0,v)$. For this special choice of configuration $ABCDEF$, System $(*)$ is reduced into two distinct equations only, namely
$$
\frac{1}{r_1^2-2a_3U+U^2}+\frac{1}{r_1^2-2a_3V+V^2}=\frac{1}{r_1^2-2a_3u+u^2}+\frac{1}{r_1^2-2a_3v+v^2},
$$
$$
\frac{1}{r_2^2-2d_3U+U^2}+\frac{1}{r_2^2-2d_3V+V^2}=\frac{1}{r_2^2-2d_3u+u^2}+\frac{1}{r_2^2-2d_3v+v^2}.
$$
For $i=1,2$, let $$\kappa_i=\frac{1}{r_i^2-2a_3u+u^2}+\frac{1}{r_i^2-2a_3v+v^2}.$$ The above equations can be rewritten, and the following System $(***)$ is obtained.
\begin{equation}
\label{ketto}
\begin{array}{llll}
(r_1^2+v^2-2a_3V)+(r_1^2+U^2-2a_3U)-\kappa_1(r_1^2+V^2-2a_3V)(r_1^2+U^2-2a_3U)=0,\\
(r_2^2+v^2-2d_3V)+(r_2^2+u^2-2d_3U)-\kappa_2(r_2^2+V^2-2d_3V)(r_2^2+U^2-2d_3U)=0.
\end{array}
\end{equation}
Writing both equations as a polynomial in $U$ gives $\omega_i(U)=\alpha_i U^2+\beta_i U + \gamma_i=0$ for $i=1,2$ where
\begin{equation}
\label{harom}
\begin{cases}
\alpha_1=\kappa_1(2Va_3-(r_1^2+V^2))+1,\\
\beta_1=2a_3(\kappa_1(V^2-2Va_3+r_1^2)-1),\\
\gamma_1=\kappa_1r_1^2(2a_3V-(r_1^2+V^2))+V^2-2Va_3+2r_1^2,\\
\end{cases}
\end{equation}
and
\begin{equation}
\label{haromplus}
\begin{cases}
\alpha_2=\kappa_2(2Vd_3-(r_2^2+V^2))+1,\\
\beta_2=2d_3(\kappa_2(V^2-2Vd_3+r_2^2)-1),\\
\gamma_3=\kappa_2r_2^2(2d_3V-(r_2^2+V^2))+V^2-2Vd_3+2r_1^2,\\
\end{cases}
\end{equation}
Since $\alpha_2\omega_1(U)-\alpha_1\omega_2(U)$ is linear in $U$ and it is equal to zero, $U$ can be written as a rational function
$$r(V)=\frac{\gamma_2\alpha_1-\gamma_1\alpha_2}{\beta_1\alpha_2-\beta_2\alpha_1}$$
with real coefficients. Thus, $V\in \mathbb{R}$ implies $U\in \mathbb{R}$.

Now, replacing $U$ by $r(V)$ in the first equation (\ref{ketto}) gives a rational function whose numerator is a polynomial $w(V)$ of degree $8$.   A straightforward computation shows that the coefficient of $V^8$ in $w(V)$ is
$$
\begin{array}{llll}
\kappa_1\kappa_2(a_3^2(4r_2^2\kappa_1\kappa_2-4\kappa_1)+d_3^2(4r_1^2\kappa_1\kappa_2-4\kappa_2)+
4a_3d_3(-\kappa_1\kappa_2(r_1^2+r_2^2)+\kappa_1+\kappa_2)+
((\kappa_1\kappa_2(r_1^2-r_2^2)+\kappa_1-\kappa_2)^2.
\end{array}
$$
It turns out that System $(***)$ has in general eight solutions, say $(u_i,v_i)$ with $i=1,\ldots 8$, where $v_1,\ldots,v_8$ are the roots of $w(V)$. Specializing $a_3,d_3,r_1,r_2,u,v$, the number of real solutions of $w(V)$ can be computed with elementary methods from the function $W=w(V)$ graph. Here we limit ourselves to the following choices: $$r_1=2,\,r_2=3,\,a_3=1,\,d_3=-1,\,u=3,\,v=-2.$$
Then $\kappa_1=19/84$ and $\kappa_2=11/72$.
Moreover, $ABC$, $DEF$ are equilateral triangles with $$A=(1,0,1), B=(\frac{1}{4} (\sqrt{6}+\sqrt{2}),\frac{1}{4} (\sqrt{6}-\sqrt{2}),1),\,\, C=(\frac{1}{4} (\sqrt{6}-\sqrt{2}), \frac{1}{4}(\sqrt{6}+\sqrt{2}),1)$$   and
 $$D=(-\sqrt{2},0,-1), E=(\frac{1}{2}(\sqrt{3}+1),\frac{1}{2}(\sqrt{3}-1),-1),\,\, F=(\frac{1}{2}(\sqrt{3}-1),\frac{1}{2}(\sqrt{3}+1),-1).$$
 A straightforward computation shows that up to a c non-zero coefficient, $\omega(V)$ is equal to
 $$(v-3)(v+2) (119089v^6 + 119089v^5 - 2707455v^4 + 1929347v^3 - 4743934v^2 - 11426976v - 3429216).$$
The third factor is a degree $6$ polynomial with two real roots, $v_1=4.1362$ and $v_2=-5.7036$, approximately. Therefore, System \ref{ketto} has two real solutions $(v_1,v_2)$ and $(v_2,v_1)$ other than the (trivial) ones $(3,-2)$ and $(-2,3)$. Thus, System $(***)$ (and hence the corresponding System $(*)$) has four real solutions.

\section{Unique Identifiability}
Some tests made with {\bf{Macaualy2}} show that assumptions (i) and (ii) are not sufficient for uniqueness identifiability. It is natural to ask whether this can be ensured by
adding a seventh equation to System $(*)$. We will show that the answer is affirmatively provided that the seven fixed points $A,B,C,D,E,F,G$ together with  $X^*,Y^*$  are chosen in $\mathbb{R}^3$
generally.

Therefore, we are led to investigate the unique identifiability problem for the extended System $(**)$ consisting  of the following seven equations
\begin{equation}
\label{s2}
\frac{1}{\|X-T\|^2} +\frac{1}{\|Y-T\|^2}=\frac{1}{\|X^*-T\|^2} +\frac{1}{\|Y^*-T\|^2} ,\quad T\in\{A,B,\ldots,F,G\}.
\end{equation}
According to Theorem~\ref{main}, System $(**)$ has finitely many solutions, and let $(X_i^*,Y_i^*)$ with $i=1,\ldots, N$ be those different from $(X^*,Y^*)$ and $(Y^*,X^*)$. Observe that the conjugate $(\overline{X_i^*},\overline{Y_i^*})$ of $(X_i^*,Y_i^*)$ is also a solution. Also, let
$X_i^*=(x_{i1}^*,x_{i2}^*,x_{i3}^*), Y_i^*=(y_{i1}^*,y_{i2}^*,y_{i3}^*)$. 

Consider the following system of equations in the unknowns $X,Y,Z$ with $X=(x_1,x_2,x_3),Y=(y_1,y_2,y_3),$ $Z=(z_1,z_2,z_3)$.
\vspace{0.2cm}
\begin{equation}
    \label{s1}
\begin{cases} \vspace{0.2cm}
{\mbox{$\frac{1}{\|X-Z\|^2} +\frac{1}{\|Y-Z\|^2}=\frac{1}{\|X^*-Z\|^2} +\frac{1}{\|Y^*-Z\|^2}$,}} \\
\vspace{0.2cm}
\prod_{i=1}^N (x_1-x_{i1}^*)\prod_{i=1}^N (x_1-y_{i1}^*)=0,
\vspace{0.2cm}
\prod_{i=1}^N (x_2-x_{i2}^*)\prod_{i=1}^N (x_2-y_{i2}^*)=0,\\
\vspace{0.2cm}
\prod_{i=1}^N (x_3-x_{i3}^*)\prod_{i=1}^N (x_3-y_{i3}^*)=0,
\vspace{0.2cm}
\prod_{i=1}^N (y_1-y_{i1}^*)\prod_{i=1}^N (y_1-x_{i1}^*)=0,\\
\vspace{0.2cm}
\prod_{i=1}^N (y_2-y_{i2}^*)\prod_{i=1}^N (y_2-x_{i2}^*))=0,
\vspace{0.2cm}
\prod_{i=1}^N (y_3-y_{i3}^*)\prod_{i=1}^N (y_3-x_{i3}^*)=0.\\
\end{cases}
\end{equation}
Note that if $(X_i^*,Y_i^*)$ is a solution of $(\ref{s2})$ for some $i$ with $1\leq i \leq N$, then
both $(X,Y,Z)=(X_i^*,Y_i^*,G)$ and $(X,Y,Z)=(Y_i^*,X_i^*,G)$ are solutions of $(\ref{s1})$. Thus,  we first investigate System $(**)$. Again, algebraic geometry provides a useful tool for our study.
\subsection{Basic facts on algebraic varieties and their projections}\label{bfag}
In the $n$-dimensional complex projective space $\mathbf{P}^n(\mathbb{C})$, $n\ge 3$, equipped with homogeneous coordinates $(\xi_1,\ldots,\xi_n,\tau)$, let $E$ be a $d$-dimensional subspace which is the intersection of the hyperplane at infinity $\tau=0$ and the $n-d-1$ hyperplanes of equations $\xi_i=0$, with $i=n-d+1,\ldots,n$. Also, let $F$ be the ($n-d-1$)-dimensional subspace of equation $\xi_1=0,\ldots,\xi_{n-d}=0$. Then $E\cap F=\emptyset$, and the projection $\pi$ with center $E$ into $F\cong \mathbf{P}^{n-d-1}(\mathbb{C})$ maps any point $P\in \mathbf{P}^n(\mathbb{C})\setminus E$ into the point $P'\in F$ uniquely determined by the elementary geometric facts that through $P$ and $E$ there passes a unique $d+1$ dimensional subspace  and that this subspace intersects $F$ in a unique point, named $P'$, and called the projection of $P$. If $V$ is a projective variety disjoint from $E$, then $\pi: V\mapsto F$ is a regular map, and the image $\pi(V)$ is closed. For more details, the reader is referred to~\cite[Chapter I]{shafaverich2013basic}.

\subsection{The unique identifiability problem for System $(**)$}
Let $V$ be the (non-empty, possible reducible) affine variety defined by the above seven equations in (\ref{s1}). To investigate $V$, our essential tool will be a projection of varieties, according to the definition recalled in subsection~\ref{bfag}. Since the image of an affine variety under a regular map need not be an affine variety, but the image of a projective variety is always a projective variety, see~\cite[Theorem 2, Chapter 1 S5]{shafaverich2013basic}; it will be useful to work with the projective closure $\overline{V}$ of $V$ defined by the seven following homogeneous equations in the $9$-dimensional complex projective space
equipped with homogeneous coordinates $(x_1:x_2:x_3:y_1:y_2:y_3:z_1:z_2:z_3:t)$.
\begin{equation}
    \label{s1h}
\begin{cases} \vspace{0.2cm}
(\|X-Z\|^2+\|Y-Z\|^2)\|tX^*-Z\|^2\|tY^*-Z\|^2=(\|tX^*-Z\|^2+\|tY^*-Z\|^2)\|X-Z\|^2\|Y-Z\|^2,\\
\vspace{0.2cm}
\prod_{i=1}^N (x_1-tx_{i1}^*)\prod_{i=1}^N (x_1-ty_{i1}^*)=0,
\vspace{0.2cm}
\prod_{i=1}^N (x_2-tx_{i2}^*)\prod_{i=1}^N (x_2-ty_{i2}^*)=0,\\
\vspace{0.2cm}
\prod_{i=1}^N (x_3-tx_{i3}^*)\prod_{i=1}^N (x_3-ty_{i3}^*)=0,
\vspace{0.2cm}
\prod_{i=1}^N (y_1-ty_{i1}^*)\prod_{i=1}^N (y_1-tx_{i1}^*)=0,\\
\vspace{0.2cm}
\prod_{i=1}^N (y_2-ty_{i2}^*)\prod_{i=1}^N (y_2-tx_{i2}^*))=0,
\vspace{0.2cm}
\prod_{i=1}^N (y_3-ty_{i3}^*)\prod_{i=1}^N (y_3-tx_{i3}^*)=0.\\
\end{cases}
\end{equation}
From the above definition, $\overline{V}$ is a real variety.  Let
$$\pi_6:
\begin{cases} {\overline{V}\rightarrow \mathbf{P}^3(\mathbb{C})},\\
 (x_1:x_2:x_3:y_1:y_2:y_3:z_1:z_2:z_3:t) \mapsto (z_1:z_2:z_3:t)
\end{cases}
$$
be the projection of $\overline{V}$ whose vertex is the $5$-dimensional subspace $E$ in the $9$-dimensional complex projective space
defined by $z_1=z_2=z_3=t=0$.
\begin{lem}
$E\cap\overline{V}= \emptyset$, that is, $\pi_6$ is well defined.
\end{lem}
\begin{proof}
Assume on the contrary that $E\cap\overline{V}\neq  \emptyset$, and let $$v=(x_{v1}:x_{v2}:x_{v3}:y_{v1}:y_{v2}:y_{v3}:z_{v1}:z_{v2}:z_{v3}:t_v)\in E\cap\overline{V}.$$
As $v\in E$ we have $t_v=0$. Therefore evaluating $(\ref{s1h})$ in $v$ yields $x_{vi}^{2N}=y_{vi}^{2N}=0$ for $i=1,2,3$. Also $z_{v1}=z_{v2}=z_{v3}=0$ by $v\in E$.
Thus, $v$ cannot define a point in $\mathbf{P}^9(\mathbb{C})$, contradicting the assumption.
\end{proof} Therefore $\pi_6(\overline{V})$ is a real projective variety of  $\mathbf{P}^3(\mathbb{C}).$
\begin{lem}
\label{lem3feb} $\dim(\pi_6(\overline{V}))<3$.
 \end{lem}
\begin{proof}
Assume on the contrary $\dim(\pi_6(\overline{V}))=3$.
Since $\pi_6(\overline{V})$ is a projective variety of $\mathbf{P}^3(\mathbb{C})$, \cite[Theorem 1.19]{shafaverich2013basic} shows that $\pi_6(\overline{V})=\mathbf{P}^3(\mathbb{C})$.
In particular, $\pi_6(\overline{V})$ is an irreducible variety of dimension $3$. Therefore, any point of $\pi_6(\overline{V})$ has a non-empty fiber.

However, we will show that there exists $\overline{ Q}$ such that $\pi_6^{-1}(\overline{Q})=\emptyset$.
For this purpose, fix two complex numbers $z_1^*,z_2^*$, and define the variable point $\overline{ Q}=\overline{Q}(z)=(z_1^*,z_2^*,z,1)$ with $z\in \mathbb{C}$. We show the existence of a finite
subset $\Gamma$ of $\mathbb{C}$ such that if $z_3\not\in \Gamma$ then $\pi_6^{-1}(\overline{Q})=\emptyset$.

Take any affine point $(X_0:Y_0:z_1^*:z_2^*:z_3:1)\in \pi_6^{-1}(\overline{Q})$ with $z_3\in \mathbb{C}$. By ($\ref{s1h}$), we have a finite set $\Sigma$ such that $(X_0,Y_0)\in \Sigma$. Let
$Q(z)=(z_1^*,z_2^*,z)$ and
\begin{equation}
\label{eqA22apr}
\begin{array}{llll}
p_3(z)&=&\|X_0-Q(z)\|^2\|Y_0-Q(z)\|^2(\|X^*-Q(z)\|^2+\|Y^*-
Q(z)\|^2)-\\
& &\|X^*-Q(z)\|^2\|Y^*-Q(z)\|^2(\|X_0-Q(z)\|^2+\|Y_0-Q(z)\|^2).
\end{array}
\end{equation}
The first equation in $(\ref{s1h})$ yields $p(z_3)=0$.
Now, the goal is to show that $z_1^*$ and $z_2^*$ can be chosen in such a way that $p(z)$ is not
 the zero polynomial for any $(X_0,Y_0)$. For such an ordered pair $(z_1^*,z_2^*)$,
 if $\Gamma$ is the set of all roots of $p(z)$ where $(X_0,Y_0)$ range over $\Sigma$
 then for every
 $z_3^*\not\in \Gamma$ we have $\pi_6^{-1}(\overline{Q})=\emptyset$ with $\overline{Q}=(z_1^*:z_2^*:z_3^*:1)$.

We need an appropriate coordinate system in $\mathbb{R}^3$ to avoid tedious computation. Let $M$ be the middle point of the segment $X^*Y^*$. In terms of vectors, $M=\ha(X^*+Y^*)$.
The translation that takes $M$ to the origin is an Euclidean transformation, so we may assume $M$ is the origin. Then $M=O$, and hence $Y^*=-X^*$.
Let $d$ be the distance between $X^*$ and $Y^*$. Then both $X^*$ and $Y^*$ lie on the sphere $\cC$ of center $O$ and radius $r=\ha d$.
The point $R^*=(r,0,0)$ is also on $\cC$. There is an Euclidean transformation $\rho$ which takes $X^*$ to $R^*$.
To show how such a $\rho$ can be found, consider the plane $\pi$ through the points $O,X^*,R^*$. In $\pi$ there exists a rotation $\rho_0$ of center $O$ which takes
$X^*$ to $R^*$. Let $\alpha$ be the angle of $\rho_0$.
Now, $\rho$ can be taken in $\mathbb{R}^3$ as the rotation of angle $\alpha$ whose axis is the line through O and orthogonal to $\pi$.
This allows us to assume $X^*=(r,0,0),Y^*=(-r,0,0)$, that is, they are two opposite points on $\cC$. So, $x_2^*=x_3^*=y_2^*=y_3^*=0$, and $y_1^*=-x_1^*=r$.
Let $X_0=(x_{01},x_{02},x_{03})$ and $Y_0=(y_{01},y_{02},y_{03})$.

In the new coordinate system in $\mathbb{R}^3$,
\begin{equation}
\label{eqA22aprbis}
\begin{array}{llll}
p_3(z)&=&((x_{01}-z_1)^2+(x_{02}-z_2)^2+\\
&&(x_{03}-z)^2+(y_{01}-z_1)^2+(y_{02}-z_2)^2+(y_{03}-z)^2)\\
& &((r-z_1)^2+z_2^2+z^2)((r+z_1)^2+z_2^2+z^2)-\\
& & ((r-z_1)^2+z_2^2+z^2+(r+z_1)^2+z_2^2+z^2)\\
& & ((x_{01}-z_1)^2+(x_{02}-z_2)^2+(x_{03}-z)^2)\\
&&((y_{01}-z_1)^2+(y_{02}-z_2)^2+(y_{03}-z)^2).
\end{array}
\end{equation}
Now assume that $p_3(z)$ is the zero polynomial.
Since $\deg(p_3(z))\leq 5$, and  the leading term of $p_3(z)$ is $-2(x_{03}+y_{03})z^5$, we have
$$y_{03}=-x_{03}.$$
Then the coefficient of $z^4$ is
$2z_1x_{01}+2z_2x_{02}+6x_{03}^2+2r^2.$
Here the coefficients of $z_1$ and $z_2$ must vanish; otherwise, we find $z_1,z_2$ such that $p_3(z)$ is not the zero polynomial. Hence
$$ y_{01}=-x_{01}\qquad y_{02}=-x_{02}.$$
Now, the coefficient of $z^3$ in $p_3(z)$ reads $16z_1 x_{01}x_{03} +16
z_2x_{02}x_{03}$. Again, the coefficients of $z_1$ and $z_2$ must vanish. Hence,
$$ x_{01}x_{03}=0, \qquad x_{02}x_{03}=0.$$
If $x_{03} \neq 0$ then $x_{01}=x_{02}=0$ whence $6x_{03}^2+2r^2=0$, a contradiction.
If $x_{03}=0$ then the coefficient of $z^4$ of $p_3(z)$ is $-2(x_{01}^2 +x_{02}^2-r^2)$. Hence, $x_{01}^2 +x_{02}^2-r^2=0$. Then
the coefficient of $z^2$ in $p_3(z)$ can also be written
$$
\begin{array}{llll}
-2x_{01}^4 - 4x_{01}^2x_{02}^2 + 4x_{01}^2z_1^2 - 4x_{01}^2z_2^2 + 16x_{01}x_{02}z_1z_2 - 2x_{02}^4 - 4x_{02}^2z_1^2 + 4x_{02}^2z_2^2 -
    4z_1^2r^2 + 4z_2^2r^2 + 2r^4=\\
   8x_{02}(2z_1z_2x_{01}-z_1^2x_{02}+z_2^2x_{02})
\end{array}
$$
If the first factor, namely $8x_{02}$, vanishes  then $x_{02}=x_{03}=0$ and $x_{01}^2=r^2$, that is, either $X_i^*=X^*, Y_i^*=Y^*$ or $X_i^*=Y^*,Y_i^*=X^*$, a contradiction. If the second factor vanishes for all $z_1,z_2$ then $x_{01}=x_{02}=0$ whence $r=0$, a contradiction.
\end{proof}
\subsection{Proof of Theorem \ref{mainbis}}
$\pi_6(\overline{V})$ defines a real affine variety $W$ in the $3$-dimensional affine space $\mathbf{R}^3$ associated to $\mathbf{P}^3(\mathbb{R})$ with respect to the plane of equation $t=0$. By Lemma \ref{lem3feb}, $\dim(W)\le 2$, i.e. $W$ is a real affine surface. This together with the remark after (\ref{s1}) shows that if $G\not\in W$, i.e. $G$ is from the open set $\mathbf{R}^3(\mathbb{R})\setminus W$, then System $(**)$ has only two trivial solutions, namely $(X^*,Y^*)$ and $(Y^*,X^*)$.

\section{Acknowledgement}
The author initiated this study during a Postdoctoral Researcher position at Altoo University. The author would like to thank Kaie Kubjas for her helpful discussions of the paper's topic.

\bibliographystyle{plain}
\bibliography{references}

\end{document}